\RequirePackage{fix-cm}
\documentclass[smallextended]{svjour3}       
\smartqed  
\usepackage{graphicx}
\usepackage{amssymb,amsmath}
\newcommand{\bd}{\bold}
\newcommand{\lbl}{\label}
\newcommand{\be}{\begin{equation}}
\newcommand{\ee}{\end{equation}}
\newcommand{\beaa}{\begin{eqnarray*}}
\newcommand{\eeaa}{\end{eqnarray*}}
\newcommand{\bea}{\begin{eqnarray}}
\newcommand{\eea}{\end{eqnarray}}

\begin{document}

\title{Limiting behavior of largest entry of random tensor constructed by high-dimensional data\thanks{Jiang's research was supported by NSF Grant No. DMS-1406279, and Xie's research was supported by  NNSF Grant No.11401169 and KRPH Grant No. 20A110001. }
}

\titlerunning{Largest entry of random tensor}        

\author{Tiefeng Jiang         \and
        Junshan Xie* 
}

\authorrunning{T. Jiang         \and
        J. Xie} 

\institute{Tiefeng Jiang \at
              School of Statistics, University of Minnesota, 224 Church Street, S. E., MN 55455, USA\\
              \email{jiang040@umn.edu}           
           \and
           Junshan Xie* \at
             School of Mathematics and Statistics, Henan University, Kaifeng, 475000, China\\
              \email{junshan@henu.edu.cn}
}

\date{Received: date / Accepted: date}

\maketitle

\begin{abstract}
Let ${X}_{k}=(x_{k1}, \cdots, x_{kp})', k=1,\cdots,n$,  be a random sample of size $n$ coming from a $p$-dimensional population. For a fixed integer $m\geq 2$, consider a hypercubic random tensor $\mathbf{{T}}$ of $m$-th order and  rank $n$ with
\begin{eqnarray*}
 \mathbf{{T}}= \sum_{k=1}^{n}\underbrace{{X}_{k}\otimes\cdots\otimes {X}_{k}}_{m~multiple}=\Big(\sum_{k=1}^{n} x_{ki_{1}}x_{ki_{2}}\cdots x_{ki_{m}}\Big)_{1\leq i_{1},\cdots, i_{m}\leq p}.
\end{eqnarray*}
Let $W_n$ be the largest off-diagonal  entry of $\mathbf{{T}}$. We derive the asymptotic distribution of $W_n$ under a suitable normalization for two cases. They are  the ultra-high dimension case with $p\to\infty$ and $\log p=o(n^{\beta})$  and  the high-dimension case with $p\to \infty$ and $p=O(n^{\alpha})$ where $\alpha,\beta>0$. The normalizing constant of $W_n$ depends on $m$ and the limiting distribution of $W_n$ is a Gumbel-type distribution involved with parameter $m$.
\end{abstract}

\section{Introduction}
\label{intro}

Let $p\geq 2$ be an integer and ${X}\in\mathbb{R}^{p}$ be a random vector.
The distribution of ${X}$ serves as a population distribution. Let ${X}_{k}=(x_{k1}, \cdots, x_{kp})', 1\leq k\leq n$, be a random sample of size $n$ from the population distribution generated by ${X}$, that is, ${X}, {X}_1, \cdots, {X}_n$ are independent random vectors with a common distribution. The data matrix $\mathbf{X}=(x_{ki})_{1\leq k\leq n,1\leq i\leq p}$ produces a hypercubic random tensor $\mathbf{T} \in \mathbb{R}^{p\times  \cdots \times p}$ with order $m$ and rank $n$ defined by
\begin{eqnarray}\lbl{new1}
 \mathbf{{T}}=
 \sum_{k=1}^{n}\underbrace{{X}_{k}\otimes\cdots\otimes {X}_{k}}_{m~multiple}=\Big(\sum_{k=1}^{n} x_{ki_{1}}x_{ki_{2}}\cdots x_{ki_{m}}\Big)_{1\leq i_{1},\cdots, i_{m}\leq p}.
\end{eqnarray}
Researchers have obtained some limiting properties of the tensor  defined in (\ref{new1}).  By using similar techniques to those in the random matrix theory,  Ambainis and Harrow \cite{r1} obtained a limiting property of the largest eigenvalue and the limiting spectral distribution of random tensors. Tieplova \cite{r18}  studied the limiting spectral distribution of the sample covariance matrices constructed by the random tensor data.  Lytova  \cite{r15} further considered the central limit theorem for linear spectral statistics of the sample covariance matrices constructed by the random tensor data. Shi {\it et al}. \cite{r17} applied limiting properties of the random tensors to an anomaly detection problem in the distribution networks.

In this paper, we will study the behavior of the largest off-diagonal  entry of the random tensor $\mathbf{T}$  when both $n$ and $p$ tend to infinity. Precisely, we will work on the asymptotic distribution of
\begin{eqnarray}
 W_{n}:=\max_{1\leq i_{1}<\cdots<i_{m}\leq p} \frac{1}{\sqrt{n}}\Big|\sum_{k=1}^{n}x_{ki_{1}}x_{ki_{2}}\cdots x_{ki_{m}}\Big|\label{1.0}
 \end{eqnarray}
as $n\rightarrow\infty$ and $p\rightarrow\infty$.

 For a fixed $m\geq 2$, when the entries of the data matrix $\bd{X}=(x_{ki})_{1\leq k\leq n,1\leq i\leq p}$ are i.i.d. random variables, we will  show that  the limiting distribution of
$W_{n}$ with a suitable normalization  is the Gumbel-type distribution involved with parameter $m$.  Two typical high-dimensional cases are considered: the ultra-high dimension with $p\to\infty$ and $ \log p=o(n^{\beta}) $ and the high-dimension with $p\to \infty$ and $p=O(n^{\alpha})$ where $\alpha,\beta>0$. In both cases we obtain the limiting distributions of $W_n$, which is different from the case that $m=2$.

When $m=2$, the tensor $\mathbf{T}=\bd{X}'\bd{X}$ turns out to be  the sample covariance matrix, which is a very popular statistic in the multivariate statistical analysis. The largest entry of the sample covariance matrix has been studied actively. In particular, assuming $n/p\rightarrow\gamma>0$ and $E|x_{11}|^{30+\epsilon}<\infty$ for some $\epsilon>0$,   Jiang \cite{r9} proved that
\begin{eqnarray}
W_{n}^{2}-4\log p +\log \log p   \xrightarrow{d}   W_{\infty}\label{1.1}
 \end{eqnarray}
where random variable $W_{\infty}$ has distribution function $F(z)=e^{-\frac{1}{\sqrt{8\pi}}e^{-z/2}},\, z\in \mathbb{R}.$ Here and later the notation ``$\xrightarrow{d}$" means ``converges in distribution to".

A sequence of results are then obtained to relax the moment condition that $E|x_{11}|^{30+\epsilon}<\infty$. For example, Zhou \cite{r20} showed that (\ref{1.1}) holds if
 \begin{eqnarray}\lbl{Zhou_nice}
x^{6}P(|x_{11}x_{12}|>x)\rightarrow 0\ \mbox{as}\ x\rightarrow\infty.
   \end{eqnarray}
Liu {\it et al.} \cite{r14} proved that  (\ref{1.1}) holds provided a weaker condition is valid, that is,
   \begin{eqnarray}\lbl{student_cucumber}
 n^{3}P(|x_{11}x_{12}|>\sqrt{n\log n} )\rightarrow 0\ \mbox{as}\ n\rightarrow\infty.
    \end{eqnarray}
Besides the above two results, Li and Rosalsky \cite{r12} and Li {\it et  al.} \cite{r10}\cite{r11}  further studied the moment condition for which (\ref{1.1}) is true. In a different direction, Liu {\it et al.} \cite{r14} obtained (\ref{1.1}) for the polynomial rate such that $p=O(n^{\alpha})$; Cai and Jiang \cite{r5} derived (\ref{1.1}) for the ultra-high dimensional case with $\log p =o(n^{\alpha})$ for some $\alpha>0.$ For the compressed sensing problems and testing problems  related to $W_n$, one is referred to   Cai and Jiang \cite{r5}, Cai  {\it et al}. \cite{r4}\cite{r6}, Xiao and Wu \cite{r19} and  Shao and Zhou \cite{r16}.

In this paper, we prove that $W_n$ with a suitable normalization converges to the Gumbel-type distribution for all $m\geq 2$. The normalizing constant and the limiting distribution all depend on $m$.
Throughout the paper, the symbols $\xrightarrow{p}$ and $\xrightarrow{d}$ mean convergence in probability and convergence in distribution, respectively. We will also denote $b_{n}=o(a_{n})$ if $\lim_{n\rightarrow\infty}{b_{n}}/{a_{n}}= 0$; the notation $b_{n}=O(a_{n})$ stands for that $\{|b_n/a_n|;\, n\geq 1\}$ is a bounded sequence;  $a_{n} \sim b_{n}$ if $ \lim_{n\rightarrow\infty}a_{n}/b_{n}=1$.

The rest of the paper is organized as follows. The main results of the paper  are stated in Section 2.
In Section 3, we  present and prove  some technical lemmas, and then prove the main results.

\section{Main Results}
\label{sec:1}
We assume that $p$ depends on $n$ and simply write $p$ for brevity of notation. In case of  possible confusion, we will write $p=p_n$. Recall $\bd{X}=(x_{ki})_{1\leq k\leq n,1\leq i\leq p}$ and assume
\bea\lbl{night_5}
& &  \{x_{ki};\, 1\leq k\leq n,1\leq i\leq p\}\ \mbox{are i.i.d.}\
 \mbox{random variables with }\ \ \ \ \ \ \ \ \ \ \ \ \ \ \ \ \nonumber\\
 &&   ~~~~Ex_{11}=0\ \mbox{and} \ Ex_{11}^{2}=1.
\eea
The quantity $W_n$ is defined in (\ref{1.0}) with $m\geq 2$. In the following theorems the limiting distribution is the Gumbel distribution with  distribution function
\bea\lbl{move_low}
F_{\theta}(z)=\exp\Big\{-\frac{1}{m!\sqrt{m\pi}}e^{-z/2}\Big\},\ \  z\in \mathbb{R}.
\eea

\begin{theorem}\lbl{th2}
Assume (\ref{night_5}) with $Ee^{t_{0}|x_{11}|^{\alpha}}<\infty$ for some $0< \alpha \leq 1$ and $t_{0}>0$. Let $m\geq 2$ be a fixed integer and  $\beta=\frac{\alpha}{2m-\alpha}$. If $p\to\infty$ and $\log p=o(n^{\beta}) $  as $n\rightarrow\infty$, then
$W_{n}^{2}-2m\log p +\log \log p    \xrightarrow{d}  \theta,$
where $\theta$ has distribution function $F_{\theta}(z)$  in (\ref{move_low}).
\end{theorem}

The above theorem studies the ultra-high dimensional case, that is, the dimension $p$ can be at an exponential order of the size $n$. Here the assumption of $Ee^{t_{0}|x_{11}|^{\alpha}}<\infty$ for some $0< \alpha \leq 1$ and $t_{0}>0$ is needed to derive the limiting distribution.
Next we  consider a popular high-dimensional case in the literature such that   $p=p_{n}$ is at most a polynomial power of  $n$. We then get the same limiting distribution for $W_{n}$ under a much weaker moment condition.

\begin{theorem}\lbl{th3}
Let $\alpha>0$  and $m\geq 2$ be constants such that $$E\big[|x_{11}|^{\tau_1}\log^{\tau_2} (1+|x_{11}|)\big]<\infty$$ with $\tau_1=4m\alpha+2$ and $\tau_2=2m\alpha+\frac{3}{2}$. If $p\to \infty$ and $p=O(n^{\alpha})$, then $W_{n}^{2}-2m\log p +\log \log p\xrightarrow{d}  \theta,$
 where $\theta$ has distribution function $F_{\theta}(z)$ in  (\ref{move_low}).
\end{theorem}

The above two theorems obviously imply the following.
\begin{corollary}
\lbl{co1} Assume the conditions from either Theorem \ref{th2} or Theorem \ref{th3} hold. Then,
\begin{eqnarray*}
 \frac{W_{n} }{\sqrt{\log p}}   \xrightarrow{p}  \sqrt{2m}.
 \end{eqnarray*}
\end{corollary}

As discussed earlier, the largest entry of a sample covariance matrix has been studied with the limiting distribution stated in (\ref{1.1}). In this paper we study the same problem for $m$-order random tensor, in which the setting is a more general.  We find that the normalizing constant of $W_n^2$ is $2m\log p -\log \log p$ and the corresponding limiting distribution is given in (\ref{move_low}). Both quantities indeed depend on $m$.
We now make some further comments below.

\medskip

1. Take $m=2$, both Theorems \ref{th2} and \ref{th3} state that
\begin{eqnarray*}
W_{n}^{2}-4\log p +\log \log p   \xrightarrow{d}  \theta \sim F_{\theta}(z)=e^{-\frac{1}{\sqrt{8\pi}}e^{-z/2}} ,
\end{eqnarray*}
which is consistent  with  (\ref{1.1}).

\medskip

2. Now, instead of studying $W_n$ from (\ref{1.0}),  we consider
\bea\lbl{Jiang_Fan}
 \tilde{W}_{n}:=\max_{1\leq i_{1}<\cdots<i_{m}\leq p} \frac{1}{\sqrt{n}}\sum_{k=1}^{n}x_{ki_{1}}x_{ki_{2}}\cdots x_{ki_{m}}.
\eea
Then, by using the same proofs except changing ``$|N(0, 1)|$" to  ``$N(0, 1)$" in (\ref{UofM}) and (\ref{NDSU}), Theorems \ref{th2} and \ref{th3} still hold with the limiting distribution ``$F_{\theta}(z)$" from  (\ref{move_low}) is replaced by ``$F(z)$", where
\beaa
F(z)=\exp\Big\{-\frac{1}{2m!\sqrt{m\pi}}e^{-z/2}\Big\},\ \  z\in \mathbb{R}.
\eeaa
Corollary \ref{co1} still holds without change if ``$W_n$" is replaced by ``$\tilde{W}_{n}$".

\medskip

3. Recently Fan and Jiang \cite{r8} studied the limiting behavior of $\tilde{W}_{n}$ from (\ref{Jiang_Fan}) with $m=2$ and with  $(x_{11}, \cdots, x_{1p})' \sim N(0, \bold{\Sigma})$, where  $\bold{\Sigma}_{ii}=1$ for each $i$ and $\bold{\Sigma}_{ij}\equiv\rho>0$ for all $i\ne j$. The limiting distribution of $\tilde{W}_{n}$ is the Gumbel distribution if $\rho$ is very small;  Gaussian if $\rho$ is large; and the convolution of the Gumbel and the Gaussian distributions if $\rho$ is in between. Such setting can also be extended to $\tilde{W}_{n}$ in (\ref{Jiang_Fan}) for any $m$ with a lengthy argument. We leave it as a future work.

\medskip

4. Assume that $m=2$ and that  $(x_{11}, \cdots, x_{1p})' \sim N(0, \bold{\Sigma})$, where $\bold{\Sigma}$ is a banded matrix.  Cai and  Jiang \cite{r5} studied $W_n$ from (\ref{1.0}) and applied their results to compressed sensing problems and tests of covariance structures. It will be interesting to see if similar dependent structures can be carried out for $W_n$  with $m\geq 3.$

\medskip

5. The proofs of Theorems \ref{th2} and \ref{th3} rely on the Stein-Chen Poisson approximation method and the moderate deviations. The major technicality comes from computing $\lambda$ and bounding $b_2$ appeared in Lemma \ref{le1}.  The major difference between our proofs here and those in the literature  is that the evaluation of $\lambda$ is more involved. Furthermore, we need a significant effort to investigate $b_2$. Due to the assumption $m\geq 3$ the dependent structure appearing in $b_2$ becomes more subtle; see Lemmas \ref{le3} and \ref{Lemma_hope} for details.

\medskip

6. Taking $m=2$ and $\alpha=1$ in Theorem \ref{th3}, the required moment condition in the theorem becomes $E\big[|x_{11}|^{10}\log^{5.5} (1+|x_{11}|)\big]<\infty.$ This is stronger than (\ref{Zhou_nice}) and (\ref{student_cucumber}). In fact, Lemma \ref{Lemma_hope}   requires this condition. It is might be possible to relax  this moment assumption in Theorem \ref{th3}. We also leave it as a future project.

\medskip

7. In the paper, the random tensor $\mathbf{T}$ is constructed by the sample of a single multivariate population. In fact, the results of Theorems \ref{th2} and \ref{th3}
can also be  extended  to the tensor constructed by the samples of several independent populations with the same dimension $p$. Let ${X}^{(l)}\in\mathbb{R}^{p};\, l=1,2,\cdots, m$ be $m$ random vectors, and the $p$ entries of ${X}^{(l)}$ be i.i.d. random variables for each $l$.   The probability  distribution of each vector generates a population distribution. For each $1\leq l \leq m$, let $(x_{k1}^{(l)}, \cdots, x_{kp}^{(l)})',\,  k=1,\cdots,n$, be  a random sample of size $n$ from the population ${X}^{(l)}$.   We then have a data matrix $\bd{X}^{(l)}=(x_{ki}^{(l)})_{1\leq k\leq n,1\leq i\leq p}$ and we  define a special hypercubic random tensor $\mathbf{T}' \in \mathbb{R}^{p\times  \cdots \times p}$ with order $m$ and rank $n$ by
\begin{eqnarray*}
 \mathbf{{T}'}=
 \Big(\sum_{k=1}^{n} x_{ki_{1}}^{(1)}x_{ki_{2}}^{(2)}\cdots x_{ki_{m}}^{(m)}\Big)_{1\leq i_{1},\cdots, i_{m}\leq p}.
\end{eqnarray*}
Denote the largest element  of $\mathbf{{T}'}$ by
\begin{eqnarray*}
 W_{n}'=\max_{1\leq i_{1}<\cdots<i_{m}\leq p} \frac{1}{\sqrt{n}}\Big|\sum_{k=1}^{n}x_{ki_{1}}^{(1)}x_{ki_{2}}^{(2)}\cdots x_{ki_{m}}^{(m)}\Big|.
 \end{eqnarray*}
 By the same argument as those in the proofs of Theorems \ref{th2} and \ref{th3}, the two theorems still hold if ``$W_n$" is replaced by ``$W_{n}'$"   under some uniform moment conditions on $x_{11}^{(l)},\, 1\leq l \leq m $.

\section{ Proofs}
\label{sec:1}
\subsection{ Some Technical Lemmas}
\label{sec:2}
We start with  some technical lemmas  useful for our proofs.  The first one is a classical Stein-Chen Poisson approximation lemma, which is frequently used in  studying behaviors of maxima of almost mutual independent random variables. The following result is a special case of Theorem 1 of Arratia {\it et al}. \cite{r3}.
\begin{lemma}
\lbl{le1} Let $\{\eta_{\alpha}, \alpha \in I\}$  be  random variables on an index set $I$ and $\{B_{\alpha}, \alpha\in I\}$ be a family of subsets of $I$, that is, for each $\alpha\in I$, $B_{\alpha}\subset I$. For any $t\in\mathbb{R}$,  set $\lambda=\sum_{\alpha\in I}P(\eta_{\alpha}>t)$. Then we have
\begin{eqnarray*}
\Big|P(\max_{\alpha\in I} \eta_{\alpha}\leq t)-e^{-\lambda} \Big|\leq (1 \wedge \lambda^{-1})(b_{1}+b_{2}+b_{3}),
\end{eqnarray*}
where
\begin{eqnarray*}
b_{1}&=&\sum_{\alpha\in I} \sum_{\beta \in B_{\alpha}} P(\eta_{\alpha}>t)P(\eta_{\beta}>t),\\
b_{2}&=&  \sum_{\alpha\in I} \sum_{\alpha \neq \beta \in B_{\alpha}} P(\eta_{\alpha}>t,\eta_{\beta}>t), \\
b_{3}&=&\sum_{\alpha\in I}\big|P(\eta_{\alpha}>t|\sigma(\eta_{\beta},\beta \notin B_{\alpha}))-P(\eta_{\alpha}>t)   \big|,
\end{eqnarray*}
and $\sigma(\eta_{\beta},\beta \notin B_{\alpha})$ is the $\sigma$-algebra generated by $\{\eta_{\beta},\beta \notin B_{\alpha}\}$. In particular,
if $\eta_{\alpha}$ is independent of $\{\eta_{\beta},\beta \notin B_{\alpha}\}$ for each $\alpha$, then $b_{3}$ vanishes.
\end{lemma}

Let $\Phi(x)=\frac{1}{\sqrt{2\pi}}\int_{-\infty}^xe^{-t^2/2}\,dt$ for $x\in \mathbb{R}.$ The following conclusion is about the moderate deviation of the partial sum of i.i.d. random variables (Linnik \cite{r13}).
\begin{lemma}
\lbl{le2} Suppose $\{\zeta, \zeta_{1}, \zeta_{2}, \cdots\}$ is a sequence of i.i.d. random variables with zero mean and $E\zeta_{i}^{2}=1$. Define $S_{n}=\sum_{i=1}^{n}\zeta_{i}$.

(1) If $Ee^{t_{0}|\zeta|^{\alpha}}<\infty$ for some $0<\alpha \leq 1$ and $t_{0}>0$, then
\begin{eqnarray*}
\lim_{n\rightarrow\infty} \frac{1}{x_{n}^{2}}\log P\Big(\frac{S_{n}}{\sqrt{n}}\geq x_{n}\Big) =-\frac{1}{2}
\end{eqnarray*}
for any $x_{n}\rightarrow\infty, x_{n}=o(n^{\frac{\alpha}{2(2-\alpha)}})$.

(2) If $Ee^{t_{0}|\zeta|^{\alpha}}<\infty$ for some $0<\alpha \leq \frac{1}{2}$ and $t_{0}>0$, then
\begin{eqnarray*}
\frac{P\big(\frac{S_{n}}{\sqrt{n}}\geq x\big)}{1-\Phi(x)}\rightarrow 1
\end{eqnarray*}
holds uniformly for $0\leq x \leq o(n^{\frac{\alpha}{2(2-\alpha)}})$.
\end{lemma}

 The following result is  Proposition 4.5 from Chen {\it et al.} \cite{r7}.
\begin{lemma}\lbl{Shao} Let $\eta_i,\, 1\leq i \leq n$, be independent random variables with $E\eta_i=0$ and $Ee^{h_n|\eta_i|}< \infty$ for some $h_n>0$ and $1\leq i \leq n$. Assume that $\sum_{i=1}^nE\eta_i^2=1$. Then
\beaa
\frac{P(\sum_{i=1}^n\eta_i\geq x)}{1-\Phi(x)}=1+C_n(1+x^3)\gamma e^{4x^3\gamma}
\eeaa
for all $0\leq x \leq h_n$ and $\gamma=\sum_{i=1}^nE\big(|\eta_i|^3e^{x|\eta_i|}\big)$, where $\sup_{n\geq 1}|C_n|\leq C$ and $C$ is an absolute constant.
\end{lemma}
\begin{proposition}\lbl{what_you} Let $\{\xi_i; i\geq 1\}$ be i.i.d. random variables with $E\xi_1=0$, $E(\xi_1^2)=1$ and $E(|\xi_1|^r)<\infty$ for some $r>2.$ Let $\{ c_n>0;\, n\geq 1\}$ be  constants with $\sup_{n\geq 1}c_n<\infty$.
Assume
\bea\lbl{Then_rain}
P\big(|\xi_1|>\sqrt{n/\log n}\,\big)=o\Big(\frac{1}{n^{1+(c_n^2/2)}\sqrt{\log n}}\Big)
\eea
as $n\to\infty$. Then, $P(S_n\geq c_n\sqrt{n\log n}\,)\sim 1-\Phi(c_n\sqrt{\log n}\,)$.
\end{proposition}

Amosova \cite{r2} derived a similar result to Proposition \ref{what_you} for independent but not necessarily identically distributed random variables. If $\{\xi_i; i\geq 1\}$ are i.i.d. random variables and  $\lim_{n\to\infty}c_n=c$, then Amosova concluded that  $P(S_n\geq c_n\sqrt{n\log n}\,)\sim 1-\Phi(c_n\sqrt{\log n}\,)$ under the condition $E(|\xi_1|^{c^2+2+\epsilon})<\infty$ for some $\epsilon>0.$ This moment condition implies (\ref{Then_rain}) by the Markov inequality
and hence our proposition holds. In particular, taking $c_n\equiv c>0$, then (\ref{Then_rain}) holds if  $E\big[|\xi_1|^{c^2+2}\log^{(c^2+3)/2}(1+|\xi_1|)\big]<\infty$.  In conclusion, for the i.i.d.  case Proposition \ref{what_you} relaxes the condition required by Amosova.

\medskip

\noindent\textbf{Proof of Proposition \ref{what_you}}. By the standard central limit theorem, as $n\to\infty$,
\beaa
\sup_{a\leq x \leq b}\Big|\frac{P(S_n\geq \sqrt{n}x)}{1-\Phi(x)}-1\Big|\to 0
\eeaa
for any real numbers $b>a$. So, without loss of generality, we will prove the conclusion under the extra assumption
\bea\lbl{melon_English}
c_n\sqrt{\log n} \to \infty
\eea
as $n\to\infty.$ The proof is divided into a few steps.

{\it Step 1: truncation}. Define $a_1=1$ and $a_n=\sqrt{n/\log n}$ for $n\geq 2.$ Denote $K=E(|\xi_1|^r).$ Set
\beaa
\xi_i' & = & \xi_iI(|\xi_i|\leq a_n)- E\big[\xi_iI(|\xi_i|\leq a_n)\big];\\
\xi_i'' & = & \xi_iI(|\xi_i|> a_n)- E\big[\xi_iI(|\xi_i|> a_n)\big]
\eeaa
for $1\leq i \leq n$. Trivially, $\xi_i=\xi_i'+\xi_i''$, $|\xi_i'|\leq 2a_n$ and
\beaa
\big|E\big[\xi_iI(|\xi_i|> a_n)\big]\big| \leq \frac{K}{a_n^{r-1}}.
\eeaa
Furthermore,
\bea\lbl{train_apple}
\mbox{Var}(\xi_i'')\leq E[\xi_i^2I(|\xi_i|> a_n)]\leq \frac{K}{a_n^{r-2}}.
\eea
 Now
\bea\lbl{people_1}
\mbox{Var}(\xi_i)=\mbox{Var}(\xi_i')+\mbox{Var}(\xi_i'')+2\mbox{Cov}(\xi_i', \xi_i'').
\eea
Use the formula that $\mbox{Cov}(U-EU, V-EV)=E(UV)-(EU)EV$ for any random variables $U$ and $V$ to see
\beaa
\mbox{Cov}(\xi_i', \xi_i'')
&=& -E\big[\xi_iI(|\xi_i|\leq a_n)\big]\cdot E\big[\xi_iI(|\xi_i|> a_n)\big]\\
& = & \big(E\big[\xi_iI(|\xi_i|> a_n)\big]\big)^2\\
& \leq & E\big[\xi_i^2I(|\xi_i|> a_n)\big]\\
& \leq & \frac{K}{a_n^{r-2}}
\eeaa
by the assumption $E\xi_i=0$,  the Cauchy-Schwarz inequality and (\ref{train_apple}). This together with  (\ref{train_apple}) and (\ref{people_1}) implies that
\bea\lbl{hen}
\mbox{Var}(\xi_1)\geq \mbox{Var}(\xi_i')\geq \mbox{Var}(\xi_1) - \frac{3 K}{a_n^{r-2}}.
\eea
Set $S_n=\sum_{i=1}^n \xi_i$, $S_n'=\sum_{i=1}^n \xi_i'$ and $S_n''=\sum_{i=1}^n \xi_i''$. Then $S_n=S_n'+ S_n''.$ Thus,
\bea\lbl{for+}
P(S_n>u)\leq P(S_n'>u-v) + P(|S_n''|>v)
\eea
for any $u>v>0$. Moreover,  $S_n'\leq S_n +|S_n''|$, we see
\beaa
P(S_n'>u+v) \leq P(S_n>u) + P(|S_n''|>v)
\eeaa
for any $u>0$ and $v>0$. This leads to
\bea\lbl{for-}
P(S_n>u)\geq P(S_n'>u+v) -P(|S_n''|>v).
\eea
From the definition of $\xi_i'$, it is easy to see that $\sup_{n\geq 1}E|\xi_1'|^{r}\leq 2^rK$.
Note that
\beaa
nE\big[|\xi_1|I(|\xi_1|> a_n)\big]\leq nE(|\xi_1|^{r})\cdot \frac{1}{a_n^{r-1}} = \frac{nK}{a_n^{r-1}}.
\eeaa
Hence
\bea\lbl{Teat}
P(|S_n''|>v) & \leq & P\Big(|\sum_{i=1}^n\xi_iI(|\xi_i|>a_n)|>v-\frac{nK}{a_n^{r-1}}\Big)\nonumber\\
& \leq & nP(|\xi_1|>a_n)
\eea
provided $v>\frac{nK}{a_n^{r-1}}$.

{\it Step 2: the tail for $S_n'$}.  Set $\sigma'^2=\mbox{Var}(\xi_1').$ Trivially, $\sigma'\to 1$ as $n\to\infty.$  Take  $\eta_i=\frac{\xi_i'}{\sqrt{n}\sigma'}$. Then $|\eta_i|\leq \frac{2}{\sigma'\sqrt{\log n}}$. Therefore we see from Lemma \ref{Shao} that
\bea\lbl{China_lid}
\frac{P(S_n'\geq \sqrt{n}\sigma'x)}{1-\Phi(x)}=1+C_n(1+x^3)\gamma e^{4x^3\gamma}
\eea
where  $\sup_{n\geq 1}|C_n|\leq C$ and $C$ is an absolute constant, and
\beaa
\gamma&\leq &\frac{n}{n^{3/2}\sigma'^{3}}\cdot E\big(|\xi'_1|^3e^{2a_nx/(\sqrt{n}\sigma')}\big)\\
& \leq & \frac{1}{\sqrt{n}\sigma'^{3}}\cdot e^{2a_nx/(\sqrt{n}\sigma')}\cdot E(|\xi'_1|^3).
\eeaa
Use the fact $\sup_{n\geq 1}E|\xi_1'|^{r}\leq 2^rK$ to see that $\sup_{n\geq 1}E(|\xi'_1|^3)\leq 2^rK$ if $r\geq 3$ by the H\"{o}lder inequality. If $2<r<3$, then write
\beaa
E(|\xi'_1|^3)&=&E(|\xi'_1|^{r}\cdot |\xi'_1|^{3-r})\\
&\leq & 2^{3-r}\Big(\frac{n}{\log n}\Big)^{(3-r)/2} \cdot E(|\xi'_1|^{r})\\
& \leq & 8K\cdot \Big(\frac{n}{\log n}\Big)^{(3-r)/2}
\eeaa
by the facts that $|\xi_1'|\leq 2a_n=2\sqrt{n/\log n}$ and that $\sup_{n\geq 1}E|\xi_1'|^{r}\leq 2^rK$. In summary, if $x=O(\sqrt{n}/a_n)$ then
\beaa
\gamma \leq
\begin{cases}
O\Big(\frac{(\log n)^{(r-3)/2}}{n^{(r/2)-1}}\Big),\ \text{if $2<r< 3$;}\\
O(n^{-1/2}),\ \text{if $r\geq 3$}
\end{cases}
\eeaa
as $n\to\infty.$ In particular, noting $a_n=\sqrt{n/\log n}$,  we know  that $\gamma \to 0$ and $x^3\gamma \to 0$ since $x=O(\sqrt{n}/a_n)$.   Consequently, we have from (\ref{China_lid}) that
\bea\lbl{leap_forward}
P(S_n'\geq \sqrt{n}\sigma'x) \sim 1-\Phi(x)
\eea
under the assumption  $x=O(\sqrt{n}/a_n).$

{\it Step 3: the tail for $S_n$}. Take $u=c_n\sqrt{n\log n}$
and $v=\frac{2nK}{a_n^{r-1}}$. Then $v/u\to 0$ as $n\to\infty.$
We still write $u$ and $v$ next sometimes for short notation. By (\ref{for+}), (\ref{for-}) and (\ref{Teat}), we have
\bea\lbl{cloth_wash1}
~~~~ P(S_n>c_n\sqrt{n\log n})\leq  P(S_n'>u-v) + nP(|\xi_1|>\sqrt{n/\log n})
\eea
and
\bea\lbl{cloth_wash2}
 ~~~~P(S_n>c_n\sqrt{n\log n})
\geq  P(S_n'>u+v) - nP(|\xi_1|>\sqrt{n/\log n}).
\eea
In what follows, we will show both $P(S_n'>u+v)$ and $P(S_n'>u-v)$ are close to $P(S_n'>u).$ Since the two arguments have no difference, we will consider them simultaneously and write $u\pm v$ for the case $u+v$ and $u-v$, respectively.     Noticing $u\pm v \sim c_n\sqrt{n\log n}$. Take $x=(u\pm v)/(\sqrt{n}\sigma')$ in (\ref{leap_forward}). Then $x=O(\sqrt{n}/a_n)$ by the assumption $\sup_{n\geq 1}c_n<\infty$. From (\ref{melon_English}), $x\sim c_n\sqrt{\log n}\to \infty.$ It follows that
\beaa
P(S_n'>u\pm v)
 & \sim  & P\Big(N(0, 1)>\frac{u\pm v}{\sqrt{n}\sigma'}\Big)\\
  &\sim & \frac{1}{c_n\sqrt{2\pi\log n}}\exp\big\{-\frac{(u\pm v)^2}{2n\sigma'^2}\big\}
\eeaa
as $n\to \infty$, where we use the fact $P(N(0,1)>x)\sim \frac{1}{\sqrt{2\pi}\,x}e^{-x^2/2}$ as $x \to\infty.$ We claim
\bea\lbl{bear_professor}
\exp\Big\{-\frac{(u\pm v)^2}{2n\sigma'^2}\Big\}\cdot \exp\Big\{\frac{c_n^2\log n}{2}\Big\}\to 1
\eea
as $n\to \infty.$ In fact, write $(u\pm v)^2=u^2 + v^2 \pm 2uv$. Then
\beaa
& & -\frac{(u\pm v)^2}{2n\sigma'^2} + \frac{c_n^2\log n}{2}\nonumber\\
& = & -\frac{u^2}{2n\sigma'^2}+ \frac{c_n^2\log n}{2} -\frac{v^2\pm 2uv}{2n\sigma'^2}\nonumber\\
& = & \frac{c_n^2(\log n)}{2}\cdot \frac{\sigma'^2-1}{\sigma'^2} + O\Big(\frac{uv}{n}\Big).
\eeaa
The assertion (\ref{hen}) says that $\sigma'^2\to 1$ and $0\leq 1-\sigma'^2\leq \frac{3K}{a_n^{r-2}}=O\big(\big(\frac{\log n}{n}\big)^{(r/2)-1}\big)$. Also,
\beaa
\frac{uv}{n}=O\Big(\frac{\sqrt{n\log n}}{a_n^{r-1}}\Big)=O\Big(\frac{(\log n)^{r/2}}{n^{(r/2)-1}}\Big).
\eeaa
It follows that
\beaa
-\frac{(u\pm v)^2}{2n\sigma'^2} + \frac{c_n^2\log n}{2}=O\Big(\frac{(\log n)^{r/2}}{n^{(r/2)-1}}\Big).
\eeaa
We then confirms (\ref{bear_professor}). Therefore,
\beaa
P(S_n'>u\pm v)\sim \frac{1}{c_n\sqrt{2\pi\log n}}\cdot \frac{1}{n^{c_n^2/2}}.
\eeaa
By the given condition,
\beaa
nP(|\xi_1|>\sqrt{n/\log n})=o\Big(\frac{1}{n^{c_n^2/2}\sqrt{\log n}}\Big).
\eeaa
Comparing these with (\ref{cloth_wash1}) and (\ref{cloth_wash2}), we arrive at \beaa
P(S_n>c_n\sqrt{n\log n}) \sim 1-\Phi(c_n\sqrt{\log n})
\eeaa
as $n\to\infty.$ \hfill $\Box$

\subsection{ Main Proofs}
\label{sec:2}
 For $1\leq s \leq m-1$, define
 \begin{eqnarray}
\xi_{k}^{(s)}=\prod_{t=1}^{s}x_{kt},~\eta_{k}^{(s)}=\prod_{t=s+1}^{m}x_{kt}, ~\zeta_{k}^{(s)}=\prod_{t=m+1}^{2m-s}x_{kt}. \label{3.2}
 \end{eqnarray}
For a number $a>0$ and a sequence of positive numbers $\{a_{n}\}$ with
 $\lim_{n\to\infty}a_{n}= a$, we define
\begin{eqnarray}
\Psi_{n}^{(s)}(a_{n})=P\Big(\big|\sum_{k=1}^{n} \xi_{k}^{(s)} \eta_{k}^{(s)}  \big|\geq a_{n}\sqrt{n\log p},\, \big|\sum_{k=1}^{n}\xi_{k}^{(s)} \zeta_{k}^{(s)} \big|\geq a_{n}\sqrt{ n\log p}  \Big)\nonumber\\\label{3.3}
\end{eqnarray}
for any $1\leq s \leq m-1$. The next is a result on $\Psi_{n}^{(s)}(a_{n})$, which is a key  step in the application of the Stein-Chen Poisson approximation to prove Theorem \ref{th2}.

\begin{lemma}\lbl{le3}  Let $\{a_{n};\, n \geq 1\}$ be a sequence of positive numbers with $\lim_{n\to\infty}a_{n}= a>0$.
Under the assumptions of  Theorem \ref{th2}, we have that $\max_{1\leq s \leq m-1}\Psi_{n}^{(s)}(a_{n})= o(p^{-a^{2}+\epsilon})$ for any $\epsilon>0.$
\end{lemma}
\textbf{Proof}. Let $u, v$ and $w>0$ be three numbers. It is easy to check that either $|u+v|>2w$ or $|u-v|>2w$ if $|u|\geq w$ and $|v|\geq w$. It then follows from (\ref{3.3}) that
\begin{eqnarray}
\Psi_{n}^{(s)}(a_{n})&\leq& P\Big(\big|\sum_{k=1}^{n} \xi_{k}^{(s)}\big( \eta_{k}^{(s)}+\zeta_{k}^{(s)} \big) \big|\geq 2a_{n}\sqrt{n\log p}\Big)\nonumber\\
&&+ P\Big(\big|\sum_{k=1}^{n}\xi_{k}^{(s)} \big( \eta_{k}^{(s)}-\zeta_{k}^{(s)} \big) \big|\geq 2a_{n}\sqrt{ n\log p}  \Big)\label{6.1}\\
&:=&A_{n}+B_{n}.\nonumber
\end{eqnarray}
For the term $A_{n}$, trivially,
\begin{eqnarray}\lbl{Tea_red}
E[\xi_{k}^{(s)}( \eta_{k}^{(s)}+\zeta_{k}^{(s)} )]=0,~
E[\xi_{k}^{(s)}( \eta_{k}^{(s)}+\zeta_{k}^{(s)} )]^{2}=2.
\end{eqnarray}
It is elementary that
\begin{eqnarray*}
  \prod_{t=1}^{m}|a_{t}|^{\alpha/m}\leq \frac{1}{m}\sum_{t=1}^{m}|a_{t}|^{\alpha}
  \end{eqnarray*}
for all $a_{t}\geq 0\, (t=1,\cdots, m)$.  Thus,  we  get
\begin{eqnarray*}
 E \exp\Big\{{t_{0}| \xi_{1}^{(s)} \eta_{1}^{(s)}  |^{\alpha/m}}\Big\} &=&
 E \exp\Big\{t_{0} \prod_{t=1}^{m}|x_{1t}  |^{\alpha/m}\Big\}\\
 &\leq & E\exp\Big\{\frac{t_{0}}{m}\sum_{t=1}^{m}|x_{1t}  |^{\alpha}\Big\}\\
 &=& \prod_{t=1}^{m}E \exp\Big\{\frac{t_{0}}{m}|x_{11}|^{\alpha}\Big\}.
\end{eqnarray*}
By assumption, $Ee^{t_{0}|x_{11}|^{\alpha}}<\infty$ for some $0< \alpha \leq 1$ and $t_{0}>0$, we see that
\begin{eqnarray}
E \exp\Big\{{t_{0}| \xi_{1}^{(s)} \eta_{1}^{(s)}  |^{\alpha/m}}\Big\}< \infty.\label{6.2}
\end{eqnarray}
Noticing $0<\alpha/m <1,$ we have
\begin{eqnarray*}
& & E \exp\Big\{\frac{1}{2}t_{0}\big| \xi_{1}^{(s)} \big( \eta_{1}^{(s)}+\zeta_{1}^{(s)} \big)  \big|^{\alpha/m}\Big\}\\
& \leq & E \Big[\exp\Big\{\frac{1}{2}t_{0}\big| \xi_{1}^{(s)} \eta_{1}^{(s)}|^{\alpha/m}\Big\}\cdot  \exp\Big\{\frac{1}{2}t_{0}\big| \xi_{1}^{(s)} \zeta_{1}^{(s)}|^{\alpha/m}\Big\}\Big]\\
& \leq & \Big[E \exp\Big\{t_{0}\big| \xi_{1}^{(s)} \eta_{1}^{(s)}|^{\alpha/m}\Big\}\Big]^{1/2}\cdot \Big[E \exp\Big\{t_{0}\big| \xi_{1}^{(s)} \zeta_{1}^{(s)}|^{\alpha/m}\Big\}\Big]^{1/2}\\
& < & \infty
\end{eqnarray*}
by the Cauchy-Schwarz inequality again. From the notation $\beta=\frac{\alpha}{2m-\alpha}$ in statement of Theorem \ref{th2}, we see    $\frac{1}{2}\cdot \frac{\alpha/m}{2-(\alpha/m)}=\frac{\beta}{2}$. An assumption implies that $a_{n}\sqrt{2\log p}=o(n^{\beta/2})$. It is easy to see that $\{\xi_{k}^{(s)}\big( \eta_{k}^{(s)}+\zeta_{k}^{(s)} \big);\, 1\leq k \leq n\}$ are i.i.d. random variables. By Lemma \ref{le2} (1) and (\ref{Tea_red}), we get that, for any sufficient small $\delta>0$,
\begin{eqnarray*}
A_{n}&\leq& P\Big(\big|\frac{1}{\sqrt{2n}}\sum_{k=1}^{n} \xi_{k}^{(s)}\big( \eta_{k}^{(s)}+\zeta_{k}^{(s)} \big) \big|\geq a_{n}\sqrt{2\log p}\Big)\\
&\leq& 2\exp \{-(1-\delta)a_{n}^{2} \log p \}\\
&=& 2p^{(\delta-1) a_{n}^{2}}.
\end{eqnarray*}
 Since $a_{n}\rightarrow a$, the above implies that, for any $\epsilon>0$,  we have
\begin{eqnarray}
A_{n} =o (p^{- a^{2}+\epsilon}) \label{6.3}
\end{eqnarray}
as $n\to \infty$. Similarly,
 \begin{eqnarray}
B_{n} =o (p^{- a^{2}+\epsilon}).\label{6.4}
\end{eqnarray}
Combining  (\ref{6.1}),  (\ref{6.3}) and (\ref{6.4}),  we complete the proof.\hfill $\Box$

\medskip

\noindent\textbf {Proof of Theorem \ref{th2}.} The asymptotic distribution of $W_{n}$  will be derived by the  Stein-Chen Poisson approximation method introduced in Lemma \ref{le1}. To do so, set $\mathbb{Z}$ be the set of integers and  $I=\{(i_{1},\cdots,i_{m})\in \mathbb{Z}^p: 1\leq i_{1}<\cdots<i_{m}\leq p\}$. For each  $\alpha=(i_{1},\cdots,i_{m})\in I$, define
\begin{eqnarray}\lbl{what_water1}
X_{\alpha}=\frac{1}{\sqrt{n}}\big|\sum_{k=1}^{n}x_{ki_{1}}x_{ki_{2}}\cdots x_{ki_{m}}\big|
\end{eqnarray}
 and
\begin{eqnarray*}
B_{\alpha}=\big\{(j_{1},\cdots,j_{m})\in I;\, \{j_{1},\cdots,j_{m}\}\cap\{i_{1},\cdots,i_{m}\}\neq \emptyset ~\mbox{but}~ (j_{1},\cdots,j_{m})\neq \alpha\big\}.
\end{eqnarray*}
Obviously, $X_{\alpha}$ is independent of $\{X_{\beta};\, \beta\in I\backslash X_{\alpha}\}$. It is easy to verify that
\bea\lbl{fix_tea}
|I|=\binom{p}{m}\ \ \mbox{and}\ \ |B_{\alpha}| \leq m^2p^{m-1}
\eea
for each $\alpha \in I.$ For  any  $z\in \mathbb{R}$,  write
 \begin{eqnarray}\lbl{what_water2}
 \nu_{p}=\Big[\log p-\frac{1}{2m}\big(\log\log p +2\log (m!\sqrt{m\pi})-z\big)\Big]^{1/2}.
 \end{eqnarray}
Notice $v_p$ may not make sense for small values of $p$. Since $p=p_n\to \infty$ as $n\to\infty$, without loss of generality, assume $v_p>0$ for all $n\geq 1$.   Set $\alpha_0=\{1,2,\cdots, m\}\in I$. By Lemma \ref{le1},
\begin{eqnarray}
\Big|P\big(\max_{\alpha \in I } X_{\alpha}\leq \sqrt{2m}\nu_{p}\big)-e^{-\lambda_{p}} \Big|\leq b_{1}+b_{2},\label{3.5}
\end{eqnarray}
where $b_1$ and $b_2$ are as in Lemma \ref{le1} and
\begin{eqnarray}
\lambda_{p}&=&
\binom{p}{m}P\big(X_{\alpha_0}>\sqrt{2m}\nu_{p}\big) \nonumber\\
&\sim &\frac{p^{m}}{m!}P\Big(\big|\sum_{k=1}^{n}x_{k1}x_{k2}\cdots x_{km}\big|>{\sqrt{2mn}\nu_{p}}\Big). \lbl{rice_egg}
\end{eqnarray}
First,  write   $\psi_{k}=x_{k1}x_{k2}\cdots x_{km},\, 1\leq k \leq n$. Then $E\psi_{k}=0$ and $E\psi_{k}^{2}=1$. The assertion (\ref{6.2}) says that $E e^{t_{0}|\psi_{1}|^{\alpha/m}}< \infty.$
Note that $\frac{\alpha}{m}\leq \frac{1}{2}$ since $0<\alpha \leq 1$ and $m\geq 2$. Moreover, $\sqrt{2m}\nu_{p}=O(\sqrt{\log p})=o(n^{\beta/2})$. By the definition of $\beta$, we know $\frac{\beta}{2}=\frac{1}{2}\cdot \frac{\alpha/m}{2-(\alpha/m)}$. Therefore it follows from   Lemma \ref{le2}(2) that
\bea
P\Big(\big|\sum_{k=1}^{n}x_{k1}x_{k2}\cdots x_{km}\big|>{\sqrt{2mn}\nu_{p}}\Big)
 &\sim & P\big(|N(0, 1)|>{\sqrt{2m}\nu_{p}}\big)\nonumber\\
 & \sim & \frac{2}{\sqrt{4\pi m}\,v_p}\cdot e^{-mv_p^2}, \lbl{UofM}
\eea
where the fact $P(N(0, 1)>x)\sim \frac{1}{\sqrt{2\pi\,}\, x}\cdot e^{-x^2/2}$ as $x\to \infty$ is used. It is easy to check that
\beaa
& & v_p\sim \sqrt{\log p};\\
& & -mv_p^2=-\log (p^m)+\frac{1}{2}\big(\log\log p +2\log (m!\sqrt{m\pi})\big)-\frac{1}{2}z
\eeaa
as $n\to\infty$. Therefore,
\begin{eqnarray}
\lambda_{p}&\sim &{\frac{p^{m}}{m!}}\frac{2}{\sqrt{4m\pi\log p}}\cdot e^{-mv_p^2} =e^{- z/2}. \label{3.6}
\end{eqnarray}
In particular, this implies that
\bea\lbl{shine_hua}
P\big(X_{\alpha}>{\sqrt{2m}\nu_{p}}\big)\sim \frac{m!}{p^m}e^{-z/2}
\eea
as $n\to\infty$ for any $\alpha\in I.$ Consequently, we have from (\ref{fix_tea}) that
\begin{eqnarray}
 b_{1} & \leq & |I|\cdot |B_{\alpha}|\cdot P\big(X_{\alpha}>\sqrt{2m}\nu_{p}\big)^{2} =O\Big(\frac{1}{p}\Big). \label{3.7}
\end{eqnarray}

Now we estimate $b_2$. First,
\bea
b_2&=&\sum_{\alpha \in I}\sum_{\beta \in B_{\alpha}}P\big(X_{\alpha}>\sqrt{2m}\nu_{p},\, X_{\beta}>\sqrt{2m}\nu_{p}\big) \nonumber\\
& = & \sum_{\alpha \in I}\sum_{s=1}^{m-1}\sum_{\beta \in I: |\beta \cap \alpha|=s}P\big(X_{\alpha}>\sqrt{2m}\nu_{p},\, X_{\beta}>\sqrt{2m}\nu_{p}\big).\lbl{Steel_19}
\eea
For $1\leq s \leq m-1$, we know
\beaa
& & \{\beta \in I: |\beta \cap \{1,2,\cdots, m\}|=s\}\\
& = & \{(i_{1},\cdots,i_{m})\in \mathbb{Z}^p: 1\leq i_{1}<\cdots<i_{m}\leq p,\, |(i_{1},\cdots,i_{m})\cap \{1,2,\cdots, m\}|=s \}.
\eeaa
Hence, $|\{\beta \in I: |\beta \cap \{1,2,\cdots, m\}|=s\}|=\binom{m}{s}\cdot \binom{p-m}{m-s}\leq m^sp^{m-s}.$ Review the notations in (\ref{3.2}) and (\ref{3.3}). In particular,
\beaa
\Psi_{n}^{(s)}(a_{n})=P\Big(\big|\sum_{k=1}^{n} \xi_{k}^{(s)} \eta_{k}^{(s)}  \big|\geq a_{n}\sqrt{n\log p},\, \big|\sum_{k=1}^{n}\xi_{k}^{(s)} \zeta_{k}^{(s)} \big|\geq a_{n}\sqrt{ n\log p}  \Big).
\eeaa
Since $x_{ij}$'s are  i.i.d. random variables, we  see that
\begin{eqnarray}\lbl{our_water}
b_{2} \leq  |I|\cdot  \sum_{s=1}^{m-1}m^sp^{m-s}\cdot \Psi_{n}^{(s)}(a_{n})
\end{eqnarray}
where
\begin{eqnarray*}
a_{n}:=\frac{\sqrt{2m}\nu_{p}}{\sqrt{\log p}}\rightarrow \sqrt{2m}
\end{eqnarray*}
as $n\rightarrow\infty$. By Lemma  \ref{le3}, for any $\epsilon>0$, we have
 \begin{eqnarray*}
\Psi_{n}^{(s)}(a_{n})\leq p^{-2m+\epsilon}
\end{eqnarray*}
as $n$ is large enough. This implies that
 \begin{eqnarray}
b_{2}&\leq & p^m\cdot m^mp^{m-1}\cdot p^{-2m+\epsilon}\nonumber\\
&=& m^mp^{-1+\epsilon}
\rightarrow 0 \label{3.8}
 \end{eqnarray}
as $n\to\infty$ for all $\epsilon \in (0, 1).$
Combining (\ref{3.5}), (\ref{3.6}), (\ref{3.7}) and (\ref{3.8}), we complete the proof. \hfill $\Box$

\medskip
The following two lemmas are prepared for the proof of Theorem \ref{th3}.
\begin{lemma}\lbl{victory_11} Let $x_{ij}$'s be as in Theorem \ref{th3} and  $\nu_{p}$ be as in (\ref{what_water2}). Define  $c_n={\sqrt{2m/\log n}\nu_{p}}$ and $\xi_1=x_{11}x_{12}\cdots x_{1m}$. Then
\bea\lbl{fan_book}
n^{1+(c_n^2/2)}\sqrt{\log n}\cdot P\big(|\xi_1|>\sqrt{n/\log n}\,\big)\to 0
\eea
as $n\to\infty$.
\end{lemma}
\noindent\textbf{Proof}. Recall $\tau_2=2m\alpha+\frac{3}{2}$ and  $\tau_1=4m\alpha +2$ and $g(x)=x^{\tau_1}\log^{\tau_2} (1+x)$ for $x\geq 0.$ Observe that
\beaa
Eg(|\xi_1|)=E\big[|x_{11}|^{\tau_1}\cdots |x_{1m}|^{\tau_1}\log^{\tau_2}(1+|x_{11}|\cdots |x_{1m}|)\big].
\eeaa
Use the inequality $1+|x_{11}|\cdots |x_{1m}|\leq (1+|x_{11}|)\cdots (1+|x_{1m}|)$ to see that
\beaa
\log^{\tau_2}(1+|x_{11}|\cdots |x_{1m}|)
& \leq & \Big[\sum_{j=1}^m\log (1+|x_{1j}|)\Big]^{\tau_2}\\
& \leq & m^{\tau_2-1}\sum_{j=1}^m\log^{\tau_2} (1+|x_{1j}|)
\eeaa
by the convex inequality. Obviously, the given condition $E\big[|x_{11}|^{\tau_1}\log^{\tau_2} (1+|x_{11}|)\big]<\infty$ implies that $E(|x_{11}|^{\tau_1})<\infty$. It follows that
\bea
Eg(|\xi_1|)
& \leq & m^{\tau_2-1}\sum_{j=1}^mE\big[|x_{11}|^{\tau_1}\cdots |x_{1m}|^{\tau_1}\log^{\tau_2} (1+|x_{1j}|)\big]\nonumber\\
& = & m^{\tau_2}E\big[|x_{11}|^{\tau_1}\log^{\tau_2} (1+|x_{11}|)\big]\cdot \big(E|x_{11}|^{\tau_1}\big)^{m-1} \nonumber\\
& < & \infty. \lbl{Process_red}
\eea
Therefore,
\beaa
P\big(|\xi_1|>\sqrt{n/\log n}\,\big)\leq \frac{Eg(|\xi_1|)}{g\big(\sqrt{n/\log n}\,\big)}.
\eeaa
Trivially, $\log (1+\sqrt{\frac{n}{\log n}})\geq \frac{1}{3}\log n$ as $n$ is sufficiently large. We then see that
\bea\lbl{cloud_fly_2}
~~~~~~g\big(\sqrt{n/\log n}\,\big)\geq 3^{-\tau_2}n^{\tau_1/2} (\log n)^{\tau_2-(\tau_1/2)}
=3^{-\tau_2}n^{2m\alpha+1}\sqrt{\log n}.
\eea
In summary,
\beaa
n^{1+(c_n^2/2)}\sqrt{\log n}\cdot P\big(|\xi_1|>\sqrt{n/\log n}\,\big) =O\big( n^{(c_n^2/2)-2m\alpha}\big).
\eeaa
The condition $p=O(n^{\alpha})$  implies that $\log p\leq \alpha\log n +O(1)$. Then we have from (\ref{what_water2}) that
\beaa
\frac{c_n^2}{2}=\frac{mv_p^2}{\log n} &\leq & \frac{1}{\log n}\cdot \Big(m\log p-\frac{1}{3}\log \log p\Big)\\
& \leq & m\alpha + \frac{1}{\log n}\cdot \Big(O(1)-\frac{1}{3}\log \log p\Big)
\eeaa
as $n$ is sufficiently large. Hence
\beaa
n^{(c_n^2/2)-2m\alpha}\leq n^{-m\alpha}\cdot \exp\Big(O(1)-\frac{1}{3}\log \log p\Big)=O\Big(\frac{1}{n^{m\alpha}(\log p)^{1/3}}\Big)
\eeaa
as $n\to\infty$. The assertion (\ref{fan_book}) is yielded. \hfill $\Box$

\begin{lemma}\lbl{Lemma_hope}
Let the assumptions of  Theorem \ref{th3} hold. Recall $\nu_{p}$ as in (\ref{what_water2}). Set $a_n=\sqrt{2m/\log p}\,\nu_p.$
Let $\Psi_{n}^{(s)}(a_{n})$ be as in (\ref{3.3}).
Then  $\max_{1\leq s \leq m-1}\Psi_{n}^{(s)}(a_{n})= O(p^{-2m+\delta})$ for any $\delta>0.$
\end{lemma}
\textbf{Proof}. It is enough to show $\Psi_{n}^{(s)}(a_{n})= O(p^{-2m+\delta})$ for each $1\leq s \leq m-1$, where $\delta>0$ is given.  Similar to (\ref{6.1}) we have that
\begin{eqnarray}
\Psi_{n}^{(s)}(a_{n})&\leq& P\Big(\big|\sum_{k=1}^{n} \xi_{k}^{(s)}\big( \eta_{k}^{(s)}+\zeta_{k}^{(s)} \big) \big|\geq 2a_{n}\sqrt{n\log p}\Big)\nonumber\\
&& \ \ \ + P\Big(\big|\sum_{k=1}^{n}\xi_{k}^{(s)} \big( \eta_{k}^{(s)}-\zeta_{k}^{(s)} \big) \big|\geq 2a_{n}\sqrt{ n\log p}  \Big)\nonumber\\
&:=&A_{n}+B_{n},\label{6.1_10}
\end{eqnarray}
where $\xi_{k}^{(s)}, \eta_{k}^{(s)}$ and $\zeta_{k}^{(s)}$ are as in (\ref{3.2}). Define $V_k=\xi_{k}^{(s)}\big( \eta_{k}^{(s)}+\zeta_{k}^{(s)} \big)/\sqrt{2}$ for $1\leq k \leq n.$ Then $\sum_{k=1}^{n}\xi_{k}^{(s)} \big( \eta_{k}^{(s)}+\zeta_{k}^{(s)} \big)=\sqrt{2}\sum_{k=1}^{n}V_k$. Observe that $V_k$'s are  i.i.d. random variables with
\begin{eqnarray}\lbl{Tea_red_11}
EV_{1}=0\ \ \mbox{and}\ \
EV_{1}^{2}=1.
\end{eqnarray}
Review $\tau_2=2m\alpha+\frac{3}{2}$, $\tau_1=4m\alpha +2$ and $g(x)=x^{\tau_1}\log^{\tau_2} (1+x)$ for $x\geq 0$ as in Lemma \ref{victory_11}. We claim $g(x)$ is a convex function on $[0, \infty).$ In fact,
\beaa
g'(x)=\tau_1 x^{\tau_1-1}\log^{\tau_2} (1+x) + \frac{\tau_2 x^{\tau_1}\log^{\tau_2-1} (1+x)}{1+x}.
\eeaa
Since $\tau_2>\frac{3}{2}$ and $\tau_1>2$, the function $\tau_1 x^{\tau_1-1}\log^{\tau_2} (1+x)$ is increasing in $x \in [0, \infty)$, and  hence its derivative is non-negative. Therefore, the convexity of $g(x)$ hinges on whether $h(x):=\frac{x^{\tau_1}\log^{\tau_2-1} (1+x)}{1+x}$ is increasing on $[0, \infty).$ Trivially,
\beaa
&& h'(x)\\
&=&\frac{1}{(1+x)^2}\Big[\underbrace{(1+x)}_{J_1}\Big(\tau_1 x^{\tau_1-1}\log^{\tau_2-1} (1+x)+ \underbrace{(\tau_2-1)\frac{x^{\tau_1}\log^{\tau_2-2} (1+x)}{1+x}}_{J_2}\Big)\\
& & \ \ \ \ \ \ \ \ \ \ \ \ \ \ \  -x^{\tau_1}\log^{\tau_2-1} (1+x)\Big]\\
& \geq & \frac{(\tau_1-1)x^{\tau_1}\log^{\tau_2-1} (1+x)}{(1+x)^2}\\
&\geq & 0
\eeaa
by using the fact $J_1>x$ and $J_2\geq 0$. Thus, $g(x)$ is convex on $[0, \infty).$ Now, by the convex property,
\beaa
Eg(|V_1|) & \leq & Eg\Big(\frac{2|\xi_{1}^{(s)} \eta_{1}^{(s)}|+2|\xi_{1}^{(s)}\zeta_{1}^{(s)}|}{2}\Big)\\
& \leq & \frac{1}{2}\big[Eg\big(2|\xi_{1}^{(s)} \eta_{1}^{(s)}|\big)+Eg\big(2|\xi_{1}^{(s)}\zeta_{1}^{(s)}|\big)\Big]\\
& = & Eg\big(2|\xi_{1}^{(s)} \eta_{1}^{(s)}|\big).
\eeaa
Since $\log (1+2x) \leq 2\log (1+x)$ for $x\geq 0$, we have $g(2x)\leq 2^{\tau_2+\tau_1}g(x)$ for $x\geq 0.$ By (\ref{3.2}), $\xi_{1}^{(s)} \eta_{1}^{(s)}=x_{11}x_{12}\cdots x_{1m}$. Consequently,
\beaa
Eg(|V_1|) \leq 2^{\tau_2+\tau_1}Eg(|x_{11}x_{12}\cdots x_{1m}|)<\infty
\eeaa
by (\ref{Process_red}). This particularly implies $E\big[g(|V_1|)I(|V_1|>\sqrt{n/\log n})\big]\to 0$.  Now,
\beaa
A_n=P\Big(\big|\sum_{k=1}^{n} V_{i}\big|\geq c_n\sqrt{n\log n}\Big)
\eeaa
where $c_n:=a_n\sqrt{2(\log p)/\log n}$. By the Markov inequality,

\beaa
P\big(|V_1|>\sqrt{n/\log n}\,\big) &\leq & \frac{E\big[g(|V_1|)I(|V_1|>\sqrt{n/\log n})\big]}{g\big(\sqrt{n/\log n}\,\big)}\\
& = & o\Big(g\big(\sqrt{n/\log n}\,\big)^{-1}\Big).
\eeaa
From (\ref{cloud_fly_2}),
\beaa
g\big(\sqrt{n/\log n}\,\big)\geq
3^{-\tau_2}n^{2m\alpha+1}\sqrt{\log n}.
\eeaa
By definition, $\lim_{n\to\infty}a_n=\sqrt{2m}$ and $a_n\leq \sqrt{2m}$ as $n$ is sufficiently large. Then
\beaa
c_n\leq 2\sqrt{m}\cdot \Big(\frac{\log p}{\log n}\Big)^{1/2}=2\sqrt{m\alpha}\Big[1+O\Big(\frac{1}{\log n}\Big)\Big]
\eeaa
by the assumption $p=O(n^{\alpha})$. This implies that
\beaa
n^{1+(c_n^2/2)}\sqrt{\log n}\cdot P\big(|V_1|>\sqrt{n/\log n}\,\big) &= & o\big( n^{(c_n^2/2)-2m\alpha}\big)\\
&= & o\big( n^{O(1/\log n)}\big)\to 0
\eeaa
 as $n\to\infty$. Therefore, it is seen from (\ref{Tea_red_11}) and Proposition \ref{what_you} that
\beaa
P\Big(\big|\sum_{k=1}^{n} V_{i}\big|\geq c_n\sqrt{n\log n}\Big) \sim 1-\Phi(c_n\sqrt{\log n}).
\eeaa
Noting $c_n\sqrt{\log n}=a_n\sqrt{2\log p}\sim \sqrt{4m\log p}$ and $1-\Phi(x)\sim \frac{1}{\sqrt{2\pi}\, x}\cdot e^{-x^2/2}$ as $x\to \infty$. It follows that
\beaa
P\Big(\big|\sum_{k=1}^{n} V_{i}\big|\geq c_n\sqrt{n\log n}\Big) =
O\Big(e^{-c_n^2(\log n)/2}\Big)=O(p^{-2m+\delta})
\eeaa
as $n\to\infty$ for any $\delta>0$. Therefore,  $A_n=O(p^{-2m+\delta})$. Similarly, $B_n=O(p^{-2m+\delta})$. The proof follows from (\ref{6.1_10}). \hfill$\square$

\medskip

\noindent\textbf{Proof of Theorem \ref{th3}.} Set $I=\{(i_{1},\cdots,i_{m})\in \mathbb{Z}^p: 1\leq i_{1}<\cdots<i_{m}\leq p\}$. For each  $\alpha=(i_{1},\cdots,i_{m})\in I$, define
\begin{eqnarray*}
X_{\alpha}=\frac{1}{\sqrt{n}}\big|\sum_{k=1}^{n}x_{ki_{1}}x_{ki_{2}}\cdots x_{ki_{m}}\big|
\end{eqnarray*}
 and
\begin{eqnarray*}
B_{\alpha}=\big\{(j_{1},\cdots,j_{m})\in I;\, \{j_{1},\cdots,j_{m}\}\cap\{i_{1},\cdots,i_{m}\}\neq \emptyset ~\mbox{but}~ (j_{1},\cdots,j_{m})\neq \alpha\big\}.
\end{eqnarray*}
Obviously, $X_{\alpha}$ is independent of $\{X_{\beta};\, \beta\in I\backslash X_{\alpha}\}$. Review (\ref{what_water1}) - (\ref{rice_egg}) in the proof of Theorem \ref{th2}. Set $\alpha_0=\{1,2,\cdots, m\}\in I$. It is seen from Lemma \ref{le1} that,
\begin{eqnarray}
\Big|P(\max_{\alpha \in I } X_{\alpha}\leq \sqrt{2m}\nu_{p})-e^{-\lambda_{p}} \Big|\leq b_{1}+b_{2},\label{3.5_10}
\end{eqnarray}
where $b_1$ and $b_2$ are as in Lemma \ref{le1} and
\begin{eqnarray*}
\lambda_{p}
&\sim &\frac{p^{m}}{m!}P\Big(\big|\sum_{k=1}^{n}x_{k1}x_{k2}\cdots x_{km}\big|>{\sqrt{2mn}\nu_{p}}\Big).
\end{eqnarray*}
Write ${\sqrt{2mn}\nu_{p}}=c_n\cdot \sqrt{n\log n}$. Immediately $c_n\to \sqrt{2m\alpha}$ as $n\to\infty$ by (\ref{what_water2}).  Set $\xi_i=x_{i1}x_{i2}\cdots x_{im}$ for $1\leq i \leq n$. Then $E\xi_1=0$,  $\mbox{Var}(\xi_1)=1$ and
\beaa
n^{1+(c_n^2/2)}\sqrt{\log n}\cdot P\big(|\xi_1|>\sqrt{n/\log n}\,\big)\to 0
\eeaa
as $n\to\infty$ by Lemma \ref{victory_11}. The assumption $E\big[|x_{11}|^{\tau_1}\log^{\tau_2} (1+|x_{11}|)\big]<\infty$ implies that $E|x_{11}|^{\tau_1}<\infty$, and hence $E|\xi_1|^{\tau_1}<\infty$ with $\tau_1=4m\alpha+2>2$. We then have from  Proposition \ref{what_you} that
\bea\lbl{NDSU}
P\Big(\big|\sum_{k=1}^{n}x_{k1}x_{k2}\cdots x_{km}\big|>{\sqrt{2mn}\nu_{p}}\Big)
 &\sim & P\big(|N(0, 1)|>{\sqrt{2m}\nu_{p}}\big)\nonumber\\
 & \sim & \frac{2}{\sqrt{4\pi m}\,v_p}\cdot e^{-mv_p^2}
\eea
as in (\ref{UofM}).  Hence,
\begin{eqnarray*}
\lambda_{p}&\sim &{\frac{p^{m}}{m!}}\frac{2}{\sqrt{4m\pi\log p}}\cdot e^{-mv_p^2} =e^{- z/2}.
\end{eqnarray*}
Immediately,
\beaa
P\big(X_{\alpha}>{\sqrt{2m}\nu_{p}}\big)\sim \frac{m!}{p^m}e^{-z/2}
\eeaa
as $n\to\infty$ for any $\alpha\in I.$ Similar to (\ref{shine_hua}) and (\ref{3.7}), we get
\begin{eqnarray}
 b_{1}  =O\Big(\frac{1}{p}\Big). \label{3.7_11}
\end{eqnarray}
Now we work on $b_2$. Recalling (\ref{Steel_19}) and (\ref{our_water}) we have
\beaa
b_2
& \leq &  |I|\cdot  \sum_{s=1}^{m-1}m^sp^{m-s}\cdot \Psi_{n}^{(s)}(a_{n})
\eeaa
where $a_{n}:=\sqrt{2m/\log p}\,\nu_{p}$ for $n\geq 1$
 and
\beaa
\Psi_{n}^{(s)}(a_{n})=P\Big(\big|\sum_{k=1}^{n} \xi_{k}^{(s)} \eta_{k}^{(s)}  \big|\geq a_{n}\sqrt{n\log p},\, \big|\sum_{k=1}^{n}\xi_{k}^{(s)} \zeta_{k}^{(s)} \big|\geq a_{n}\sqrt{ n\log p}  \Big)
\eeaa
and $\xi_{k}^{(s)}$, $\eta_{k}^{(s)}$ and $\zeta_{k}^{(s)}$ are as in (\ref{3.2}).  By Lemma  \ref{Lemma_hope},
 \begin{eqnarray*}
\max_{1\leq s \leq m-1}\Psi_{n}^{(s)}(a_{n})\leq O(p^{-2m+\delta})
\end{eqnarray*}
as $n\to\infty$ for any $\delta>0$. By using the fact $|I|\leq p^m$, we have
 \begin{eqnarray*}
b_{2}&\leq & p^m\cdot m^mp^{m-1}\cdot p^{-2m+\delta}\nonumber\\
&=& m^mp^{-1+\delta}
\rightarrow 0
 \end{eqnarray*}
as $n\to\infty$ for any $\delta \in (0,1)$. This joining with (\ref{3.5_10})-(\ref{3.7_11}) completes the proof.  \hfill $\Box$

\section*{Acknowledgements}
The authors would like to thank the referee and an associate editor for their constructive comments and suggestions that have led to improvements in the paper.



%
%



\end{document}